\documentclass[10pt]{article}
\usepackage[pagebackref=true,colorlinks=true, linkcolor=blue ,  citecolor=blue, urlcolor=black]{hyperref}
\usepackage[dvipsnames]{xcolor}  
\usepackage{fullpage}
\usepackage{amsmath}
\usepackage{amsthm}

\usepackage{mathrsfs}
\usepackage{fancyhdr}

\usepackage{amssymb}
\usepackage{mathtools}
\usepackage{xcolor}
\usepackage{algorithm}
\usepackage{algorithmic}
\usepackage{bbold}
\usepackage{cite}
\usepackage{latexsym}
\usepackage{booktabs}
\usepackage{bm}
\usepackage{graphicx}
\usepackage{textcomp}
\usepackage{subfigure}
\usepackage[toc,page]{appendix}
\usepackage{url}
\usepackage{hyperref}
\usepackage{multicol}
\usepackage{tikz}

\let\oldtextit\textit 
\renewcommand\emph[1]{\oldtextit{\color{RoyalBlue}#1}}

\theoremstyle{definition}
\newtheorem{definition}{Definition}[section]

\newtheorem{prop}[definition]{Proposition}

\newtheorem{theorem}[definition]{Theorem}

\theoremstyle{remark}
\newtheorem{remark}[definition]{Remark}

\makeatletter
\newcommand{\address}[1]{\gdef\@address{#1}}
\newcommand{\email}[1]{\gdef\@email{\url{#1}}}
\newcommand{\urladdr}[1]{\gdef\@urladdr{\url{#1}}}
\newcommand{\@endstuff}{\par\vspace{\baselineskip}\noindent\small
\begin{tabular}{@{}l}\scshape\@address\\\textit{E-mail address:}
\@email
\\\textit{URL: }\@urladdr\end{tabular}}
\AtEndDocument{\@endstuff}
\makeatother
\title{The NumericalCertification package in Macaulay2}
\author{Kisun Lee}
\date{}
\address{Department of Mathematics, University of California San Diego,\\ 9500 Gilman Drive, La Jolla, CA 92093}
\email {kil004@ucsd.edu}\urladdr{https://klee669.github.io}

\begin{document}

\maketitle

\begin{abstract}
    The package \texttt{NumericalCertification} implements methods for certifying numerical approximations of solutions for a given system of polynomial equations. For certifying regular solutions, the package implements Smale's $\alpha$-theory and Krawczyk method. For a singular solution, we implement soft verification using the iterative deflation method. We demonstrate the functionalities of the package focusing on interaction with current numerical solvers in \texttt{Macaulay2}.
\end{abstract}

\section{Introduction}

Systems with polynomial equations arise in many fields in mathematics and applied science. Specially, the interest on polynomial systems turns into a problem in algebraic geometry that finding all isolated solutions for a given polynomial system. Due to a recent development in \emph{numerical algebraic geometry} (e.g.\ see \cite{SommeseWampler:2005}), a family of numerical algorithms called the \emph{homotopy continuation} gains popularity as a way to find solutions for a polynomial system. There are several known  implementations \texttt{Bertini} \cite{BHSW06}, \texttt{HomotopyContinuation.jl} \cite{breiding2018homotopycontinuation}, \texttt{Hom4PS-3} \cite{chen2014hom4ps}, \texttt{NumericalAlgebraicGeometry} \cite{leykin2011numerical} and \texttt{PHCpack} \cite{verschelde1999algorithm} which are all widely used.

One remark for the homotopy continuation algorithm is that its output is not certified. It means that numerical approximations obtained by the algorithm might not satisfy the users depending on their purposes. We say that a numerical approximation is \emph{certified} if a compact region that contains a unique solution can be obtained from the given approximation by applying a sort of algorithm. For a certified approximation, the unique solution contained in the compact region is called an \emph{associated solution} and the given approximation is called an \emph{approximate solution}. We call this series of algorithms \emph{numerical certification}.

As implementations for numerical certification for a polynomial system, we point out \texttt{alphaCertified} \cite{hauenstein2011alphacertified} and a function \texttt{certify} in the software \texttt{HomotopyContinuation.jl} \cite{breiding2020certifying}. \texttt{alphaCertified} implements \emph{Smale's $\alpha$-theory} \cite{smale1986newton} as a way for numerical certification \cite{hauenstein2012algorithm}. On the other hand, \texttt{HomotopyContinuation.jl} exploits \emph{Krawczyk method} \cite{krawczyk1969newton} using \emph{interval arithmetic} \cite{Moore:1977} as a tool for certification. 

The package \texttt{NumericalCertification} in \texttt{Macaulay2}\cite{M2} executes regular solution certification using both $\alpha$-theory and Krawczyk method. A preferred method can be chosen as an option by users. As an improved version of the software presented in \cite{lee2019certifying}, it includes soft verification for a singular solution using the idea of the \emph{deflation} method \cite{leykin2006newton}. Finally, the package provides an interface to the software \texttt{alphaCertifed}. 

The rest of the paper consists of two sections. In the next section, we discuss the required preliminaries for certification. The implementation details are given in the last section.

\section{Preliminaries}

In this section, we review the concepts used for numerical certification. Smale's $\alpha$-theory and a combination of Krawczyk method and interval arithmetic are used as methods for regular solution certification. We introduce the deflation method for a notion providing an idea for singular solution certification.
\subsection{Smale's \texorpdfstring{$\alpha$}--theory}

Consider an $n\times n$-square system $F$, i.e.\ a system with $n$ polynomial equations with $n$ variables.
For a point $x\in \mathbb{C}^n$, recall the \emph{Newton operator} $N_F(x)$ defined like the following :
\[N_F(x)=\left\{\begin{array}{ll}
x-F'(x)^{-1}F(x)
     & \text{if }F'(x)\text{ is invertible}, \\
x     & \text{otherwise}
\end{array}\right.\] 
We say that a sequence $\{N_F^k(x)\}_{k=1}^\infty$ \emph{converges quadratically} to an associated solution $x^\star$ of $F$ if for every $k\in \mathbb{Z}_>$,
\[\left\|N_F^k(x)-x^\star\right\|\leq \left(\frac{1}{2}\right)^{2^k-1}\|x-x^\star\|.\]
In this case, $x$ is an approximate solution for $F$.
When $F'(x)$ is not invertible, we say $x$ is an approximate solution if and only if $F(x)=0$. The $\alpha$-theory provides a certificate for the quadratic convergence of a given point. The certificate is obtained from the three auxiliary parameters : 
\begin{equation*}
	\begin{array}{ccl}
	\alpha(F,x)& := & \beta(F,x)\gamma(F,x)\\
	\beta(F,x) & := & \|x-N_F(x)\|=\|F'(x)^{-1}F(x)\|\\
	\gamma(F,x) & := & \sup\limits_{k\geq 2}\left\|\frac{F'(x)^{-1}F^{(k)}(x)}{k!}\right\|^{\frac{1}{k-1}}
	\end{array}
\end{equation*}
where $F^{(k)}(x)$ in the definition of $\gamma(F,x)$ is a symmetric tensor whose components are the $k$-th partial derivatives of $F$, see \cite[Chapter 5]{MR783635}. The norm in $\beta(F,x)$ is the usual Euclidean norm and the norm in $\gamma(F,x)$ is the operator norm on $S^k\mathbb{C}^n$ (for details, see \cite{hauenstein2011alphacertified}).  When $F'$ is not invertible at $x$, we define $\alpha(F,x)=\beta(F,x)=\gamma(F,x)=\infty$.
The contents of $\alpha$-theory are summarized below.
\begin{theorem}[c.f.\   \cite{blum2012complexity,hauenstein2012algorithm}]
Let $F:\mathbb{C}^n\rightarrow \mathbb{C}^n$ be a square polynomial system with a point $x\in \mathbb{C}^n$. Then, 
\begin{enumerate}
    \item if $\alpha(F,x)<\frac{13-3\sqrt{17}}{4}$, then $x$ is an approximate solution for $F$, and $\|x-x^\star\|\geq 2\beta(F,x)$ where $x^\star$ is an associated solution to $x$.
    \item if $\alpha(F,x)<0.03$ and $\|x-y\|<\frac{1}{20\gamma(F,x)}$ for a point $y$, then $x$ and $y$ are both approximate solutions for $F$ to the same solution $x^\star$, and $\|x-x^\star\|\leq \frac{1}{20\gamma(F,x)}$.
    \item if $\|x-\overline{x}\|>4\beta(F,x)$ for the conjugate $\overline{x}$ of $x$, then $x^\star$ is not real.
\end{enumerate}
\end{theorem}

For an implementation of $\alpha$-theory, the step for computing (or bounding) $\gamma$ is required. 
For a degree $d$ polynomial $f=\sum_{|\nu|\leq d}a_\nu x^\nu$, we recall that the \emph{Bombieri-Weyl norm} is defined as 
\begin{equation*}
    \|f\|^2=\frac{1}{d!}\sum\limits_{|\nu|\leq d}\nu!(d-|\nu|)! |a_\nu|^2.
  \end{equation*}
For a system of polynomials $F=\{f_1,\dots,f_n\}$, we define a norm for the system
$$
\|F\|^2=\sum_{i=1}^n\|f_i\|^2.
$$
Let $d_i=\deg f_i$ for each $i=1,\dots, n$ and $d=\max d_i$.  For a point $x\in\mathbb{C}$, define $\|(1,x)\|^2=1+\sum_{i=1}^n|x_i|^2$, and we let $\Delta_F(x)$ be the diagonal matrix with entries $\Delta_F(x)_{ii}:=\sqrt{d_i} \|(1,x)\|^{d_i-1}$.
Combining all these, a bound for $\gamma(F,x)$ is given as follows:

\begin{prop}\cite[Proposition 5]{hauenstein2012algorithm}
Let $F$ be a square system of polynomials and $x\in \mathbb{C}^n$ be a point. Suppose that $F'(x)$ is nonsingular.  Define
\begin{equation}\label{eq:mu}
\mu(F,x):=\max\left\{1,\|F\|\|F'(x)^{-1}\Delta_F(x)\|\right\}
\end{equation}
where the norm in $\|F'(x)^{-1}\Delta_F(x)\|$ is the operator norm. Then, 
\begin{equation}\label{eq:gammaBound}
\gamma(F,x)\leq\frac{\mu(F,x)d^{\frac{3}{2}}}{2\|(1,x)\|}.
\end{equation}
\end{prop}

\subsection{Interval arithmetic and Krawczyk method}

Interval arithmetic is introducing arithmetic operators between intervals to achieve conservative results on numerical computations. For an operator $\odot$ with intervals $[a,b]$ and $[c,d]$ over $\mathbb{R}$, we define
$[a,b]\odot[c,d]=\{x\odot y\mid x\in [a,b] ,y\in [c,d]\}$.
The real interval arithmetic can be extended over $\mathbb{C}$. For describing an interval over the complex, we use two intervals to construct an interval box $[a_1,b_1]+i[a_2,b_2]$ containing numbers in $\mathbb{C}$. Then, complex interval arithmetic can be done similarly using arithmetic over the complex numbers. A set of complex intervals is denoted by $\mathbb{IC}$. Likewise, a set of $n$-dimensional complex interval boxes is denoted by $\mathbb{IC}^n$. For a function $F:\mathbb{C}^n\rightarrow \mathbb{C}^m$ and an interval box $I\in \mathbb{IC}^n$, we define an \emph{interval extension} of $F$ as a set containing the image of all $F$ on $I$, and it is denoted by $\square F(I)$.

Krawczyk method is a combination of interval arithmetic and generalized Newton's method to get certificates for the existence and uniqueness of a solution for a square system in a given interval. Suppose that a square differentiable system $F:\mathbb{C}^n\rightarrow \mathbb{C}^n$ is given with an interval $I\in \mathbb{IC}^n$. 
Let $Y$ be an $n\times n$-invertible matrix and $x$ be a point in $I$. Then, we define \emph{Krawczyk operator} centered at $x$ like the following :
\[K_{x,Y}(I):=x-YF(x)+(I_n-Y\square F'(I))(I-x)\]
where $I_n$ is the $n\times n$ identity matrix.
Then, the following theorem summarizes the propositions required for interval arithmetic-based certification.
\begin{theorem}[c.f.\  \cite{krawczyk1969newton,burr2019effective, breiding2020certifying}]\label{thm:Krawczyk}
Suppose that $F:\mathbb{C}^n\rightarrow\mathbb{C}^n$ is a square differentiable system with a given interval extension $\square F(I)$ on an interval $I$. For an $n\times n$-invertible matrix $Y$ and a point $x$,
\begin{enumerate}
\item if a root $x^\star$ of $F$ is in $I$, then $x^\star \in K_{x,Y}(I)$.
\item \label{thm:KrawczykExistence} if $K_{x,Y}\subset I$, then $I$ contains a root $x^\star$ of $F$.
\item \label{thm:KrawczykUniqueness} if $I$ contains a root of $F$ and $\sqrt{2}\|I_n-Y\square F'(I)\|<1$, then the root $x^\star$ in $I$ is unique.
\item if $I$ contains a root of $F$, $\sqrt{2}\|I_n-Y\square F'(I)\|<1$ and a set of conjugates $\{\overline{y}\mid y\in K_{x,Y}(I)\}$ for $K_{x,Y}(I)$ is contained in $I$, then the root $x^\star$ in $I$ is unique and real.
\end{enumerate}
Here, $\|I_n-Y\square F'(I)\|$ is the maximum operator norm of the interval matrix $I_n-Y\square F'(I)$ under the max-norm.
\end{theorem}
Note that the invertible matrix $Y$ is chosen for minimizing $\|I_n-Y\square F'(I)\|$. In an actual implementation, a natural choice for $Y$ can be $F'(m(I))^{-1}$ where $m(I)$ is the midpoint of the box $I$.

\begin{remark}
In general, interval arithmetic certification allows working with less precision than $\alpha$-theory. On the other hand, $\alpha$-theory shows a better convergence rate to an actual solution. An example in \cite[Section 5.1]{burr2019effective} shows a comparison between two methods.
\end{remark}

\subsection{The deflation method}

A deflation is a series of method to reinstate the quadratic convergence of Newton iteration for an isolated singular solution of a system of equations. The basic idea is introducing more equations to construct an augmented system with reduced singularity (e.g.\ multiplicity). For an isolated singular solution $x^\star$ of a square system $F=\{f_1,\dots, f_n\}$, define $\kappa = \dim\ker F'(x^\star)$. Then, for a randomly chosen vector $B=\begin{bmatrix}b_1,\dots,b_n\end{bmatrix}^\top\in\mathbb{C}^{n\times 1}$ from the kernel of $F'(x^\star)$, an augmented system 
\begin{equation}\label{eq:augmentedSystem}
G= \begin{bmatrix} F\\ F'B\end{bmatrix}\end{equation}
has a solution $x^\star$ with a lower multiplicity than that of $F$. It is known that a singular solution is regularized within finitely many iterations by applying iterative deflation \cite{leykin2006newton}.

\section{Implementation details}

The package \texttt{NumericalCertification} is designed to interplay with other numerical solvers in \texttt{Macaulay2}, for example, \texttt{NumericalAlgebraicGeometry} \cite{leykin2011numerical}, \texttt{Bertini} \cite{bates2013bertini} or \texttt{PHCpack} \cite{gross2011phcpack}. Hence, the package supports \texttt{PolySystem} and \texttt{AbstractPoint} types of input.

The most direct way to use the package is \texttt{certifySolutions}.
It takes a polynomial system and a list of numerical solutions as input.
	\begin{verbatim}
  i1 : needsPackage "NumericalCertification"
  i2 : R = CC[x1,x2,y1,y2];
  i3 : F = polySystem {3*y1 + 2*y2 -1, 3*x1 + 2*x2 -3.5, x1^2 + y1^2 -1, x2^2 + y2^2 -1};
  i4 : sols = solveSystem f; -- a list of numerical solutions
  i5 : c = certifySolutions(F,sols);	
\end{verbatim}
It supports three strategies as options, \texttt{alphaTheory}, \texttt{intervalArithmetic} and \texttt{alphaCertified}, and the default value is \texttt{alphaTheory}. The function returns \texttt{MutableHashTable} and we can peek the output to see the certification results.
The option \texttt{alphaTheory} returns the values of $\alpha(F,x)$ for each numerical solution (in an order of the input), the list of certified distinct, real, regular and singular solutions. The list of non-certified solutions also returned and it may be certified again after refinement.
\begin{verbatim}
  i6 : peek c
  o6 = MutableHashTable{alphaValues => {2.07811e-30, 1.97421e-40}
         certifiedDistinct => {{.652548, .771177, .757747, -.63662}, 
              {.95437, .318445, -.298627, .947941}}
         certifiedReal => {{.652548, .771177, .757747, -.63662}, 
              {.95437, .318445, -.298627, .947941}}
         certifiedRegular => {{.652548, .771177, .757747, -.63662}, 
              {.95437, .318445, -.298627, .947941}}
         certifiedSingular => {}
         nonCertified => {}
         }
\end{verbatim}
The option \texttt{intervalArithmetic} returns the list of Krawczyk operators for certified real, regular and singular solutions. For the list of non-certified solutions, it returns the input interval boxes used for certification.
\begin{verbatim}
  i7 : c = certifySolutions(F,sols, Strategy=>"intervalArithmetic");
  i8 : peek c
  o8 = MutableHashTable{certifiedReal => {|  [.95437,.95437] + 
              [-1.33962e-27,1.33962e-27]*ii ... }
         certifiedRegular => {|  [.95437,.95437] + 
              [-1.33962e-27,1.33962e-27]*ii ... }
         certifiedSingular => {}
         nonCertified => {}
         }
\end{verbatim}
For the option \texttt{alphaCertified}, we need to set a path to the software installed. It might require reloading the package. When we run \texttt{certifySolutions} with the option, it runs the software and creates files in the directory where \texttt{alphaCertified} is installed.
\begin{verbatim}
  i9 : loadPackage("NumericalCertification",
           Configuration=>{"ALPHACERTIFIEDexec"=>"path/to/alphaCertified/"},
           Reload=>true)	
  i10 : certifySoltions(F,sols, Strategy=>"alphaCertified")

    alphaCertified v1.3.0 (October 16, 2013)
    Jonathan D. Hauenstein and Frank Sottile
           GMP v6.2.1 & MPFR v4.1.0

                      ...
    
\end{verbatim}

What follows is the implementation details of each strategy used in the package.

\subsection{\texorpdfstring{$\alpha$}--theory certification}

In this section, we look into functions executed in certification with $\alpha$-theory. These functions can also be used separately.

The function \texttt{computeConstants} takes a polynomial system $F$ and a numerical point $x$, and computes three parameter values $\alpha(F,x),\beta(F,x)$ and $\gamma(F,x)$. For the value of $\gamma(F,x)$, the upper bound given in (\ref{eq:gammaBound}) is used instead. We use the Frobenius norm to bound the matrix norm used in $\mu(F,x)$ (\ref{eq:mu}). The polynomial system and its Jacobian are evaluated by a straight-line program (see \cite[Section 4.1]{burgisser2013algebraic}) implemented in the package \texttt{SLPexpressions} \cite{slpexpressions} for a faster evaluation. Note that the polynomial system and the point must be in the same coefficient ring.
\begin{verbatim}
  i11 : x = point{{.652548, .771177, .757747, -.63662_CC}};
  i12 : computeConstants(F,x)
  o12 = (1.16708e-10, 5.22384e-13, 223.414)
  o12 : Sequence
\end{verbatim}
The function \texttt{certifyRegularSolution} certifies the given solution by checking the inequality $\alpha(F,x)<\frac{13-3\sqrt{17}}{4}$. It returns true if input satisfies the inequality, false otherwise.
\begin{verbatim}
  i13 : certifyRegularSolution(F,x)
  o13 = true
\end{verbatim}
The function \texttt{certifyDistinctSolutions} takes a polynomial system and two points as input. It returns false if given two points converge to the same actual solution of the system, otherwise true.
\begin{verbatim}
  i14 : y = point{{.95437, .318445, -.298627, .947941_CC}};
  i15 : certifyDistinctSolutions(F,x,y)
  o15 = true
\end{verbatim}
A given solution may converge to a solution over the real numbers even though it is a complex-valued solution. The function \texttt{certifyRealSolution} checks if a given solution corresponds to a real solution or not.
\begin{verbatim}
  i16 : x = point{{.652548, .771177, .757747, -.63662+0.001*ii}};
  i17 : certifyRealSolution(F,x)
  o17 = true
\end{verbatim}
The function \texttt{alphaTheoryCertification} takes a polynomial system and a list of numerical solutions and runs all aforementioned functions at once according to the algorithm established in \cite[Section 2.2]{hauenstein2012algorithm}. Unlike \texttt{certifySolutions}, it does not execute singular solution certification.
\begin{verbatim}
  i18 : sols = {x,y};
  i19 : c = alphaTheoryCertification(F,sols);
  i20 : peek c
  o20 = MutableHashTable{alphaValues => {.000223414, 1.04693e-10}}
          certifiedDistinct => {x, y}
          certifiedReal => {x, y}
          certifiedRegular => {x, y}
\end{verbatim}
Finally, the package supports the exact arithmetic over the rational numbers or Gaussian rationals for $\alpha$-theory certification. For example, constants $\alpha(F,x),\beta(F,x)$ and $\gamma(F,x)$ can be computed over the Gaussian rationals as follows :
\begin{verbatim}
  i21 : CR = QQ[i]/ideal(i^2+1); -- a ring of Gaussian rationals
  i22 : R = CR[x1,x2,y1,y2];
  i23 : F = polySystem {3*y1 + 2*y2 -1, 3*x1 + 2*x2 -7/2, x1^2 + y1^2 -1, x2^2 + y2^2 -1};
  i24 : x = point(sub(matrix{{5/9,3/4,3/4,-1/2}},CR));
  i25 : computeConstants(F,x)
          73052652544805089   9731461   60054828392
  o25 =  (-----------------, ---------, -----------)
           8695980754208352  303595776   229146291
  o25 : Sequence
\end{verbatim}

\subsection{Interval arithmetic certification}
For interval arithmetic certification, the package provides a type of intervals over the complex numbers.
The function \texttt{intervalCCi} returns a complex interval from a pair of real intervals representing real and imaginary part respectively. 
\begin{verbatim}
  i26 : I1 = intervalCCi(interval(.8,.9),interval(-0.1,0.1))
  o26 = [.8,.9] + [-.1,.1]*ii
  o26 : CCi
\end{verbatim}
When only one real interval is given as input, it returns a complex interval with the zero interval for its imaginary part.
\begin{verbatim}
  i27 : I2 = intervalCCi(interval(.2,.3))
  o27 = [.2,.3] + [0,-0]*ii
  o27 : CCi
\end{verbatim}
The package supports a basic interval arithmetic for the complex intervals and matrices with complex interval entries.
\begin{verbatim}
  i28 : I1 + I2
  o28 = [1,1.2] + [-.1,.1]*ii
  o28 : CCi
  i29 : I1 * I2
  o29 = [.16,.27] + [-.03,.03]*ii
  o29 : CCi
  i30 : I1^3
  o30 = [.486,.756] + [-.244,.244]*ii
  o30 : CCi
  i31 : M = matrixCCi{{I1,I2},{I2,I1}}
  o31 = |  [.8,.9] + [-.1,.1]*ii [.2,.3] + [0,-0]*ii |
        |  [.2,.3] + [0,-0]*ii [.8,.9] + [-.1,.1]*ii |
  o31 : CCiMatrix
  i32 : M^2
  o32 = |  [.67,.91] + [-.18,.18]*ii [.32,.54] + [-.06,.06]*ii |
        |  [.32,.54] + [-.06,.06]*ii [.67,.91] + [-.18,.18]*ii |
  o32 : CCiMatrix
\end{verbatim}
The function \texttt{pointToInterval} construct an interval box from a given point. This function helps to make interval input without defining them separately. There are two ways to use the function. The first is inputting a point and a desired radius for an interval box. Then, it returns an interval box with the given radius centered at the given point.
\begin{verbatim}
  i33 : x = point{{-1.6,-1.3*ii}};
  i34 : I = pointToInterval(x,1)
  o34 = |  [-2.6,-.6] + [-1,1]*ii [-1,1] + [-2.3,-.3]*ii |
  o34 : CCiMatrix
\end{verbatim}
In many cases, a proper radius can be different depending on the polynomial system or the accuracy of the approximation. Running the function \texttt{pointToInterval} with a polynomial system and a point as input returns an interval box with a radius estimating the distance between the point and the convergence limit by using Newton-Kantorovich theorem.
\begin{verbatim}
  i35 : R = CC[x,y];
  i36 : F = polySystem {x^2 + y^2 -1, x - y^2};
  i37 : x = point{{-1.61803, -1.27202*ii}};
  i38 : I = pointToInterval(F,x)
  o38 = |  [-1.61803,-1.61803] + [-8.75505e-15,8.75291e-15]*ii
              [-1.80849e-14,1.80894e-14] + [-1.27202,-1.27202]*ii |
  o38 : CCiMatrix
\end{verbatim}
The function \texttt{krawczykOperator} computes Krawczyk operator from a given polynomial system and an interval box or a point. When a point is given as input, it computes Krawczyk operator from the interval obtained by \texttt{pointToInterval(F,x)}.
\begin{verbatim}
  i39 : krawczykOperator(F,I)
  o39 = |  [-1.61803,-1.61803] + [-2.2629e-25,2.2629e-25]*ii 
              [-1.65999e-25,1.65999e-25] + [-1.27202,-1.27202]*ii |
  o39 : CCiMatrix
  i40 : krawczykOperator(F,x)
  o40 = |  [-1.61803,-1.61803] + [-2.2629e-25,2.2629e-25]*ii 
              [-1.65999e-25,1.65999e-25] + [-1.27202,-1.27202]*ii |
  o40 : CCiMatrix
\end{verbatim}
The function \texttt{krawczykTest} checks if Krawczyk operator satisfies \ref{thm:KrawczykExistence} and \ref{thm:KrawczykUniqueness} of Theorem \ref{thm:Krawczyk}. It returns true if it is (hence the given interval is certified to contain a solution uniquely), false otherwise.
\begin{verbatim}
  i41 : krawczykTest(F,I)
  o41 = true
\end{verbatim}
Finally, the function \texttt{krawczykRealnessTest} certifies if a given interval corresponds to a real solution to the system or not. It returns true if the given interval contains a unique real solution to the system, false otherwise.
\begin{verbatim}
  i42 : krawczykRealnessTest(F,I)
  o42 = false
  i43 : y = point{{.618034, -.786151}}; -- a real solution
  i44 : krawczykRealnessTest(F,y)
  o44 = true
\end{verbatim}
As shown in \texttt{i44}, both functions \texttt{krawczykTest} and \texttt{krawczykRealnessTest} also take a point as input.

\subsection{Singular solution certification}
The method of iterated deflation is exploited for singular solution certification. For a polynomial system and a numerical approximation, we consider a subsystem $$\hat{F}=F+F'B=\left\{f_i+\sum_{j=1}^n\frac{\partial f_i}{\partial x_j}b_j\right\}_{i=1,\dots,n}$$
of the overdetermined system $G$ given in (\ref{eq:augmentedSystem}) obtained by the deflation. The $\alpha$-theory or interval arithmetic certification is applied on the square subsystem and the given numerical solution. If the numerical solution is still singular, then we construct an augmented system and take a square subsystem repeatedly. Since applying the deflation on a singular solution must terminate within finitely many iterations, the given numerical approximation becomes an approximation of a regular solution of the square subsystem.

It is possible to produce a false positive result as $\hat{F}$ is obtained by squaring-up the overdetermined system from a randomly chosen vector $B$. However, it can recover the quadratic convergence of Newton iteration for a singular solution with probability $1$, and so it can be used as soft verification of a singular solution.

Singular solution certification is done by the function \texttt{certifySingularSolution} to a given polynomial system and a numerical solution.
\begin{verbatim}
  i45 : F = polySystem {x^2 + y, x^3 - y^2};
  i46 : x = point{{1e-7,2e-7*ii}};
  i47 : certifySingularSolution(F,x)
  o47 = true
\end{verbatim}
Both strategies \texttt{alphaTheory} and \texttt{intervalArithmetic} are available as options. The function executes the iterated deflation until the given singular solution is regularized. Therefore, the function might not terminate if a poor approximation is given. To prevent this, if a user knows the number of iterations required in advance, it can be provided as input.
\begin{verbatim}
  i48 : certifySingularSolution(F,x,1) -- an insufficient number of iterations
  o48 = false
  i49 : certifySingularSolution(F,x,2)
  o49 = true
\end{verbatim}

\section*{Acknowledgements}
The author is grateful to Anton Leykin for encouraging the project and suggesting various improvements. The author would also like to thank Michael Burr and Thomas Yahl for help on the package, and the organizers of the Macaulay2 workshop at Cleveland State University where the project was greatly improved.

\bibliography{ref.bib}

\begin{thebibliography}{BGLR13}

\bibitem[BCS13]{burgisser2013algebraic}
Peter B{\"u}rgisser, Michael Clausen, and Mohammad~A Shokrollahi.
\newblock {\em Algebraic complexity theory}, volume 315.
\newblock Springer Science \& Business Media, 2013.

\bibitem[BCSS12]{blum2012complexity}
Lenore Blum, Felipe Cucker, Michael Shub, and Steve Smale.
\newblock {\em Complexity and real computation}.
\newblock Springer Science \& Business Media, 2012.

\bibitem[BGLR13]{bates2013bertini}
Daniel~J Bates, Elizabeth Gross, Anton Leykin, and Jose~Israel Rodriguez.
\newblock Bertini for macaulay2.
\newblock {\em arXiv preprint arXiv:1310.3297}, 2013.

\bibitem[BHSW]{BHSW06}
Daniel~J. Bates, Jonathan~D. Hauenstein, Andrew~J. Sommese, and Charles~W.
  Wampler.
\newblock Bertini: Software for numerical algebraic geometry.
\newblock Available at bertini.nd.edu with permanent doi:
  dx.doi.org/10.7274/R0H41PB5.

\bibitem[BLL19]{burr2019effective}
M.~Burr, K.~Lee, and A.~Leykin.
\newblock Effective certification of approximate solutions to systems of
  equations involving analytic functions.
\newblock In {\em Proceedings of the 2019 on International Symposium on
  Symbolic and Algebraic Computation}, pages 267--274, 2019.

\bibitem[BRT20]{breiding2020certifying}
Paul Breiding, Kemal Rose, and Sascha Timme.
\newblock Certifying zeros of polynomial systems using interval arithmetic.
\newblock {\em arXiv preprint arXiv:2011.05000}, 2020.

\bibitem[BT18]{breiding2018homotopycontinuation}
Paul Breiding and Sascha Timme.
\newblock Homotopy{C}ontinuation.jl: A package for homotopy continuation in
  {J}ulia.
\newblock In {\em International Congress on Mathematical Software}, pages
  458--465. Springer, 2018.

\bibitem[CDLS]{slpexpressions}
Justin Chen, Timothy Duff, Anton Leykin, and Mike Stillman.
\newblock {SLP}expressions.m2: a {M}acaulay2 package for straight line programs
  and expressions for evaluation circuits.

\bibitem[CLL14]{chen2014hom4ps}
Tianran Chen, Tsung-Lin Lee, and Tien-Yien Li.
\newblock {H}om4{PS}-3: a parallel numerical solver for systems of polynomial
  equations based on polyhedral homotopy continuation methods.
\newblock In {\em International Congress on Mathematical Software}, pages
  183--190. Springer, 2014.

\bibitem[GPV11]{gross2011phcpack}
Elizabeth Gross, Sonja Petrovi{\'c}, and Jan Verschelde.
\newblock {PHC}pack in {M}acaulay2.
\newblock {\em arXiv preprint arXiv:1105.4881}, 2011.

\bibitem[GS02]{M2}
D.~R. Grayson and M.~E. Stillman.
\newblock Macaulay2, a software system for research in algebraic geometry.
\newblock Available at \url{http://www.math.uiuc.edu/Macaulay2/}, 2002.

\bibitem[HS11]{hauenstein2011alphacertified}
Jonathan~D Hauenstein and Frank Sottile.
\newblock alpha{C}ertified: Software for certifying numerical solutions to
  polynomial equations.
\newblock Available at
  \url{http://math.tamu.edu/~sottile/research/stories/alphaCertified}, 2011.

\bibitem[HS12]{hauenstein2012algorithm}
J.~D. Hauenstein and F.~Sottile.
\newblock Algorithm 921: alpha{C}ertified: certifying solutions to polynomial
  systems.
\newblock {\em ACM Transactions on Mathematical Software (TOMS)}, 38(4):28,
  2012.

\bibitem[Kra69]{krawczyk1969newton}
Rudolf Krawczyk.
\newblock Newton-algorithmen zur bestimmung von nullstellen mit
  fehlerschranken.
\newblock {\em Computing}, 4(3):187--201, 1969.

\bibitem[Lan83]{MR783635}
Serge Lang.
\newblock {\em Real analysis}.
\newblock Addison-Wesley Publishing Company, Advanced Book Program, Reading,
  MA, second edition, 1983.

\bibitem[Lee19]{lee2019certifying}
Kisun Lee.
\newblock Certifying approximate solutions to polynomial systems on
  {M}acaulay2.
\newblock {\em ACM Communications in Computer Algebra}, 53(2):45--48, 2019.

\bibitem[Ley11]{leykin2011numerical}
A.~Leykin.
\newblock Numerical algebraic geometry.
\newblock {\em Journal of Software for Algebra and Geometry}, 3(1):5--10, 2011.

\bibitem[LVZ06]{leykin2006newton}
A.~Leykin, J.~Verschelde, and A.~Zhao.
\newblock Newton's method with deflation for isolated singularities of
  polynomial systems.
\newblock {\em Theoretical Computer Science}, 359(1-3):111--122, 2006.

\bibitem[Moo77]{Moore:1977}
R.~E. Moore.
\newblock A test for existence of solutions to nonlinear systems.
\newblock {\em SIAM Journal on Numerical Analysis}, 14(4):pp. 611--615, 1977.

\bibitem[Sma86]{smale1986newton}
S.~Smale.
\newblock Newton's method estimates from data at one point.
\newblock {\em The Merging of Disciplines: New Directions in Pure, Applied, and
  Computational Mathematics}, 1986.

\bibitem[SW05]{SommeseWampler:2005}
A.~Sommese and C.~Wampler.
\newblock {\em The Numerical Solution of Systems of Polynomials Arising in
  Engineering and Science}.
\newblock World Scientific, 2005.

\bibitem[Ver99]{verschelde1999algorithm}
Jan Verschelde.
\newblock Algorithm 795: {PHC}pack: A general-purpose solver for polynomial
  systems by homotopy continuation.
\newblock {\em ACM Transactions on Mathematical Software (TOMS)},
  25(2):251--276, 1999.

\end{thebibliography}
\bibliographystyle{alpha}

\end{document}